\documentclass[%
]{mpi2015-cscpreprint}

\usepackage{graphicx}
\usepackage{graphics}

\usepackage{epstopdf}
\usepackage[caption=false]{subfig}
\usepackage{amsmath}

\newcommand{\diagdots}[3][-25]{%
\rotatebox{#1}{\makebox[0pt]{\makebox[#2]{\xleaders\hbox{$\cdot$\hskip#3}\hfill\kern0pt}}}%
}
\usepackage{amsfonts}
\usepackage{amssymb}
\usepackage[english]{babel}
\usepackage{pgfplots}
\usepackage{nomencl}
\makenomenclature
\usetikzlibrary{matrix}
\usepackage{parskip}
\usepackage{mymacros}
\usepackage{amsthm}
\usepackage{tikz}
\usepackage{enumerate}
\usetikzlibrary{calc,patterns,decorations.pathmorphing,decorations.markings}

\usepackage{cite}
\usepackage{esvect} 
\usepackage{soul} 
\usepackage{tikz}
\tikzstyle{vertex}=[circle, draw, inner sep=0pt, minimum size=6pt]
\usepackage{float}
\usepackage{mathtools}
\usepackage{todonotes}
\usepackage{booktabs}
\usepackage{ifdraft}
\usepackage{algorithm}
\usepackage{algpseudocode}
\usepackage{setspace}
\usepackage{tablefootnote}
\usepackage{color}
\usepackage{tikz}
\usetikzlibrary{shapes.geometric, arrows}
\tikzstyle{startstop} = [rectangle, rounded corners, minimum width=3cm, minimum height=0.6cm,text centered, draw=black, fill=blue!20]
\tikzstyle{process} = [rectangle, rounded corners, minimum width=3cm, minimum height=0.6cm, text centered, draw=black, fill=blue!20]
\tikzstyle{errest} = [rectangle, rounded corners, minimum width=3cm, minimum height=0.6cm, text centered, draw=black, fill=blue!20]
\tikzstyle{arrow} = [thick,->,>=stealth]

\renewcommand{\hat}{\widehat}
\renewcommand{\tilde}{\widetilde}

\newcommand{\vvec}{\mathsf{v}}
\newcommand{\wvec}{\mathsf{w}}

\usepackage[square,sort&compress,comma,numbers]{natbib}
\bibpunct[, ]{[}{]}{,}{n}{,}{,}
\newcounter{mymac@matlab}
\newcommand{\matlab}{MATLAB%
\ifnum\value{mymac@matlab}<1%
\textsuperscript{\textregistered}%
\setcounter{mymac@matlab}{1}%
\fi%
}

\theoremstyle{plain}

\newtheorem{remark}{Remark}

\theoremstyle{definition}

\theoremstyle{remark}

\usepackage{amssymb}
\usepackage{amsthm}
\crefformat{figure}{Fig.~#2#1#3}
\crefformat{section}{Section~#2#1#3}
\crefformat{table}{Table~#2#1#3}
\crefformat{algorithm}{Algorithm~#2#1#3}

\begin{document}


\title{A Scalable Framework for Simultaneously Optimizing Damper Positions and Viscosities}

\author[$ \ast \ast~1~\dagger$ ]{M. Ugrica Vukojević}
\author[$\ast~2$]{P. Goyal}
\author[$\ast \ast ~3$]{Z. Tomljanović}
\affil[$\ast \ast$]{School of Applied Mathematics and Informatics, 
University Josip Juraj Strossmayer, Trg Ljudevita Gaja 6, 31000 Osijek, Croatia}
\affil[$\ast$]{Independent researcher, Heilbronn, Germany}

\affil[1]{\email{mugrica@mathos.hr}, \orcid{0000-0002-8585-8967}}
\affil[2]{\email{pkgpawan@gmail.com}, \orcid{0000-0003-3072-7780}}
\affil[3]{\email{ztomljan@mathos.hr}, \orcid{0000-0002-3239-760X}}
\affil[$\dagger$]{Most of this work was completed while this author was with
the Max Planck Institute for Dynamics of Complex Technical Systems, Magdeburg.}

\shorttitle{A scalable framework for damping optimization}
\shortauthor{M. Ugrica Vukojević, P. Goyal, Z. Tomljanović}


\abstract{
We consider damping optimization for vibrating systems described by a second-order differential equation.
The goal is to determine the damping positions and viscosity values of the dampers in the system so that they produce the lowest average total energy.
We propose a new framework for determining the number and positions of the dampers, using re-weighted $l_1$ minimization and pruning techniques.
The efficiency and performance of our new approach are verified and illustrated in several numerical examples.
}
\novelty{We aim to determine the ``optimal'' number of dampers, their positions, and corresponding viscosities by applying different techniques to different parts of the problem. More specifically, our contribution here consists of the new approach to viscosity optimization via penalized objective, where we aim to determine which dampers' positions are unimportant in our system. Additionally, after the viscosity optimization is done, we apply the pruning technique to determine the optimal number of dampers in the system.}

\maketitle

\section{Introduction}
We consider a damped linear vibrational system described by the
differential equation
\begin{equation}\label{eq:system}
    M \ddot{q} + C(\vvec) \dot{q} + K q = 0 \,, \qquad
    q(0)=q_0, \quad \dot{q}(0)=\dot{q}_0 \,,
\end{equation}
where $M,C(\vvec),K \in \mathbb{R}^{n \times n}$ (called mass, damping, stiffness matrix,
respectively) are real, symmetric matrices with $M, K$ being
positive definite and $C(\vvec)=C_\intern+C_\extern(\vvec)$. Additionally, $C_\intern$ is positive definite and represents the internal damping. It is often modelled as a small multiple of Rayleigh or proportional damping, i.e., $C_\intern=\alpha M +\beta K$, where $\alpha, \beta \geq 0$, or as a small multiple of the critical damping. In this paper, we focus on the latter convention of internal damping, i.e.,
\begin{equation} \label{eq:Cint}
C_\intern = 2 \alpha M^{1/2}\sqrt{M^{-1/2}KM^{-1/2}}M^{1/2},
\end{equation}
where $\alpha > 0$. For more details on critical damping, see, e.g., \cite{Adh06, Ves11, KuzTT12}. The external damping matrix $C_\extern(\vvec)$ includes the geometry of damping positions and their corresponding viscosities, i.e., 
\begin{equation}\label{eq:Cext}
C_\extern(\vvec)= \vvec_1 g_{1} g_{1}^\T+\vvec_2 g_{2} g_{2}^\T+\cdots+\vvec_s
g_{s} g_{s}^\T = G\diag{\vvec_1,\vvec_2,\dots,\vvec_s}G^\T,
\end{equation}
where the matrix $G=[g_1,g_2,\dots, g_s]$ consists of the vectors $g_{j}$, for $j=1,\ldots, s$, where $s$ is the number of dampers in the system. These vectors determine the geometry and the positions of the dampers. In the standard case $g_{j}=e_{i}$ for some $i\in \{1,\dots, n\}$, where $e_i$ is the $i$-th canonical vector, which can be interpreted as a grounded damper on the $i$-th mass. The vector $\vvec=(\vvec_1,\ldots,\vvec_s)\in \mathbb{R}^s_+$ collects the nonnegative viscosities, where $\mathbb{R}^s_+$ denotes the set of $s$-dimensional vectors with real nonnegative entries.

The problem of damping optimization has been widely studied in the last few decades, and has typically been considered in two different settings: when the system~\eqref{eq:system} is homogeneous, and when it is non-homogeneous, i.e., excited by some external force. In the most general context, the problem can be stated as follows. For given mass and stiffness matrices, we aim to determine an external damping matrix such that unwanted vibrations decay as fast as possible. For a fixed number of dampers $s$, this amounts to determining the optimal $s$ positions of the dampers together with their viscosities. For homogeneous systems a further, higher-level question arises: how should one choose the number of dampers $s$ so that a given optimization criterion is minimized up to a desired accuracy? This question can be combined with the former one of determining the positions and the respective viscosities.

Answering these questions requires a suitable optimization criterion, and the choice of criterion depends strongly on the application. An overview of different damping optimization criteria can be found in \cite{Ves11}. Criteria for the non-homogeneous case with fixed positions, where the system is additionally excited, have also been proposed in \cite{Ves11}. Recently, in \cite{Tom23}, the author considered the scenario of a periodic excitation force, while in \cite{Kuztt16}, the authors investigated the energy over an arbitrary time interval. For multiple-input multiple-output (MIMO) systems with fixed positions, damping optimization can also be performed by minimizing standard system norms, such as the $\cH_2$ or $\cH_\infty$ norm. For more details on these, we refer to \cite{morBeaGT20,morBenKTetal16,BlaCGetal12,morTomBG18,morTomV20}.

In this paper, we focus on a homogeneous case and consider minimizing the average total energy of the system, which has been considered, e.g., in \cite{morBenTT11a, morBenTT13, PTT17, CoxNRetal04, Ves11}.  In both \cite{morBenTT11a, morBenTT13}, the authors use a dimension reduction approach to accelerate optimization. Furthermore, \cite{morBenTT13} deals with the damping of selected eigenfrequencies, i.e., $\omega$ for which $Kx = \omega^2 M x$ for some $x \neq 0$, obtained from the undamped system, that is, from \eqref{eq:system} with $C(\vvec)=0$. The authors in \cite{CoxNRetal04, Ves11} give the theoretically optimal damping, which in the modal coordinates reads $\Phi^\T C(\vvec) \Phi = 2\Omega$. This choice is not feasible in practice. The main reason is that $C_\extern(\vvec)$ in modal coordinates does not have the desired structure. For instance, one would typically expect external damping modelled by a low-rank matrix; therefore, in modal coordinates this matrix cannot be diagonal and positive definite. In \cite{PTT17}, the authors propose an approximation of the average total energy for the case when we are close to a modally damped mechanical system.
The main reason is to accelerate the calculation of the average total energy and therefore accelerate viscosity optimization. On the other hand, in \cite{Pal2025}, the authors exploit the structure to efficiently compute the average total energy.

On the other hand, one can also consider an eigenvalue-based criterion to damp the resonant frequencies, as discussed in \cite{morGraMQetal16}, where the spectral abscissa criterion is minimized. The spectral abscissa is defined as 
\begin{equation}
  	\alpha_{\mathrm{MCK}}(\vvec) =\max\limits_{\lambda(\vvec) \in \Lambda(\vvec)} \mathrm{Re} \, \lambda(\vvec),
\end{equation}
where $\lambda(\vvec)$ are solutions to the following quadratic eigenvalue problem (QEP):
\begin{linenomath}
\begin{equation}
	\label{eq:qep_eq}
	(\lambda(\vvec)^2M+\lambda(\vvec) C(\vvec)+K)x(\vvec)=0,
\end{equation}
\end{linenomath}
and $\Lambda(\vvec)$ denotes the spectrum of~\eqref{eq:qep_eq}.

Another problem related to the homogeneous cases is the frequency isolation problem, investigated in  \cite{JakMTetal21, MorE16}. In \cite{MorE16}, undamped vibrational systems have been considered, i.e., $C(\vvec) = 0$, whereas in \cite{JakMTetal21}, damped systems were considered, i.e., $C(\vvec) \neq  0$, and undesirable resonant frequency bands were known \emph{a priori}.

Damping optimization is in general carried out either by choosing the damping positions, or by choosing the viscosity values (i.e., by optimizing $\vvec \in \mathbb{R}^s_+$), or by doing both simultaneously. Viscosity values are optimized in most of the papers mentioned above, e.g., in \cite{JakMTetal21, PTT17, Kuztt16, morBeaGT20}. Computing the optimal damping positions, on the other hand, is a very challenging problem, studied in \cite{GurM92, KanPTT19, PrzUVTB2026}. No efficient algorithm for it is known so far, although some heuristics can be found in, e.g., \cite{KanPTT19}.
One further approach to determining optimal damping positions is the ``brute force'' search, in which all possible damping configurations are considered and the viscosities are optimized for each configuration. This approach is computationally very demanding, because of the number of possible configurations. One more question arises here: what is the ``optimal'' number of dampers?
By the ``optimal'' number of dampers we mean the smallest number for which the relative error in the function value is not larger than a given tolerance. Since the global minimum is generally unknown, we use the best function value found during the optimization as the reference value in the computation of this relative error.

In this paper, we aim to determine the ``optimal'' number of dampers, their positions, and the corresponding viscosities by applying different techniques to different parts of the problem. More specifically, our contribution consists of a new approach to viscosity optimization via a penalized objective, by means of which we identify the damper positions that are unimportant for the system. Additionally, once the viscosity optimization is complete, we apply a pruning technique to determine the optimal number of dampers in the system.

The paper is organized as follows. In \cref{sec:setting}, we state the minimization problem in more detail, while in \cref{sec:position_opt}, we present our full algorithm for damping optimization together with some implementation details. In \cref{sec:nume_ex}, we validate our new approach and compare it with the ``brute force'' approach. \cref{sec:conclusion} concludes the paper.

\section{Problem setting} \label{sec:setting}
In this work, we utilize the average total energy as our optimization criterion. In order to write it explicitly, we need to transform the $n$-dimensional second-order system~\eqref{eq:system} into a $2n$-dimensional first-order differential equation.  To that end, first note that $M$ and $K$ are symmetric positive definite matrices; hence, there exists a matrix $\Phi$, which simultaneously diagonalizes $M$ and $K$, i.e.,
\begin{equation} \label{eq:eigdecomp}
\Phi^\T K \Phi = \Omega^2  \,, \quad  \Phi^\T M \Phi = I \,,
 \end{equation}
 where
\begin{equation} \label{eq:omege}
\Omega=\diag{\omega_1, \ldots, \omega_n} , \qquad \omega_1 \geq
\cdots \geq \omega_n\,,
\end{equation}
are the undamped frequencies (eigenfrequencies), i.e., the numbers $\omega_j$ for which $Kx = \omega_j^2 M x$ has a nontrivial solution $x$, as is also apparent from \eqref{eq:eigdecomp}; since $M$ and $K$ are positive definite, they are real and positive. Note that they are ordered \emph{decreasingly}, so that $\omega_n$ is the lowest eigenfrequency; repeated eigenfrequencies are permitted.

Furthermore, it can also be shown that $\Phi$  diagonalizes internal damping, i.e., $\Phi^\T C_\intern\Phi = 2\alpha\Omega$, see \cite{Ves11}. Thus, with substitutions $y_1(t;\vvec)= \Omega\Phi^{-1}q, \,y_2(t;\vvec)=\Phi^{-1}\dot{q}$, and $y = [y_1;y_2]$, we can transform system \eqref{eq:system} into 
\begin{equation}
\dot{y}(t;\vvec)=\bA(\vvec)\,y(t;\vvec),  
\end{equation}
where 
\begin{equation} \label{eq:matrixA}
\bA(\vvec) = \begin{bmatrix}
                 0      &  \Omega      \\[2pt]
                -\Omega & -\Phi^\T( C_\intern+C_\extern(\vvec)) \Phi
              \end{bmatrix}.
\end{equation}

In our setting, optimization is focused on optimizing the matrix $C_\extern(\vvec)$. We want to determine the optimal number of dampers and optimize positions and viscosities. For the $s$ dampers that define the external damping in \eqref{eq:Cext}, we want to determine the number of required damping positions and corresponding viscosities such that the average  total energy is minimal.

The total average energy can be calculated as $\trace{\bX(\vvec)}$, where $\bX(\vvec)$ is the solution of the following Lyapunov equation
\begin{equation}          \label{eq:total_avg_en}
   \bA(\vvec)\bX(\vvec) + \bX(\vvec)\bA(\vvec)^\T = - \bZ,
   \qquad \bZ = \begin{bmatrix} E_\cS & 0 \\ 0 & E_\cS \end{bmatrix} \,.
\end{equation}

Here $\cS \subseteq \{1,\ldots,n\}$ is the index set of those eigenfrequencies that are to be damped, and
\begin{equation} \label{eq:selector}
E_\cS = \sum_{i \in \cS} e_i e_i^\T \in \mathbb{R}^{n \times n}
\end{equation}
is the corresponding diagonal selector matrix.

Thus, $\bZ$ is a diagonal matrix that determines which of the undamped frequencies should be damped, and the matrix $\bA(\vvec)$ is given in \eqref{eq:matrixA}. If we are interested in damping all of the undamped eigenfrequencies, then $\cS = \{1,\ldots,n\}$ and $\bZ = \bI$. Moreover, it can be shown that the  average total energy is equivalent to the $\cH_2$ norm when the input-output matrices are chosen appropriately. For more details, see, e.g., \cite{NakTT19, CoxNRetal04, Ves11}.

Now, determination of the optimal damping positions and viscosities is equivalent to the minimization of
the trace of the solution of the corresponding Lyapunov equation, i.e.,
\begin{equation*}  \label{eq:traceZ}
  \min_{\vvec \in \mathbb{R}_+^s} \trace{\bX(\vvec)}
\end{equation*}

Thus, the objective function is
\begin{equation}\label{eq:min_f}
f(\vvec) = \trace{\bX(\vvec)} 
\end{equation} 
and we have $s$ variables to optimize.

Our goal is to avoid ending up with all $s$ dampers active in the system. Instead, we aim to minimize the number of dampers used, while ensuring that the objective function does not increase significantly. Ideally, the optimization process then yields $s^*$ positive viscosities, with the remaining ones equal to zero up to a prescribed tolerance, so that the optimal number of dampers is given by $s^*$.

\section{Optimizing damper positions and viscosities}\label{sec:position_opt}
In this section, we discuss a problem formulation that aims at optimizing damper positions and their corresponding viscosities simultaneously.

For ease of notation, let $\Xi=(M, K, G, \alpha, s, \mathbf{Z})$ be the tuple consisting of system matrices $M, K, G,$ coefficient $\alpha$ from internal damping, $s$ the initial dimension of vector $\vvec$, and matrix $\mathbf{Z}$ from equation \eqref{eq:total_avg_en}, which together uniquely determine our objective function $f$ given in \eqref{eq:min_f}. From now on, function $f$ will be denoted as $f_\Xi$ indicating the dependencies on $\Xi$.

\subsection{Penalized objective function}
To optimize both the number and placement of dampers, we aim to modify the total average energy criterion by introducing a penalization term. Consequently, the objective function is defined as

\begin{equation}\label{eq:optimization_l0}
	\tf_\Xi(\vvec) = f(\vvec) +  \gamma \|\vvec\|_{l_0}= \trace{\bX(\vvec)} +  \gamma \|\vvec\|_{l_0},
\end{equation}
where $\bX(\vvec)$ is the solution to the Lyapunov equation \eqref{eq:total_avg_en} with
\[
\bA(\vvec) = \begin{bmatrix}
	0      &  \Omega      \\
	-\Omega & -2\alpha\Omega-\Phi^\T C_\extern(\vvec) \Phi
\end{bmatrix},
\]
and $C_\extern(\vvec)$ is given by \eqref{eq:Cext}. Here, $\gamma$ is a fixed parameter that determines the strength of the penalization term. The $l_0$-norm, denoted by $\|\cdot\|_{l_0}$, counts the number of nonzero elements in a vector. By including $\|\vvec\|_{l_0}$ in the objective function, we can simultaneously determine both the positions and viscosities of the dampers. Specifically, for each $j \in  \{1,\dots, s\}$, a nonzero value, $\vvec_j \neq 0$, indicates the presence of a damper at geometry $g_j$ as defined in \eqref{eq:Cext} with viscosity $\vvec_j$. Conversely, $\vvec_j = 0$ implies that no damper is present at geometry $g_j$. 

The $l_0$ ``norm'' is computed as
\begin{equation} \label{eq:norm_0}
	\|\vvec\|_{l_0} = \sum\limits_{j=1}^s \II(\vvec_j \neq 0),
\end{equation}
where $\II(\cdot)$ is the indicator function defined by
\begin{equation} \label{eq:indicatior_func}
	\II(u \neq 0)=
	\begin{cases}
		1, &\text{if $u\neq 0$},\\
		0, &\text{if $u=0$}
	\end{cases}.
\end{equation}
If the problem \eqref{eq:optimization_l0} could be solved efficiently, it would yield the desired solution. However, this optimization is particularly challenging due to its non-convex nature and the combinatorial complexity arising from the presence of many local minima. Therefore, a common approach is to relax the $l_0$-norm, replacing it by a convex surrogate, most often the $l_1$-norm. Under certain conditions, the solution to this relaxed $l_1$ problem coincides with that of the original combinatorial $l_0$ problem in the linear case. We discuss the $l_1$-based relaxation, as described in \cite{CWB08}, in more detail in the next subsection. Additionally, we employ a further relaxation, using $|\vvec_j|<\texttt{tol}_0$ in place of $\vvec_j = 0$, i.e., values within this tolerance are treated as zero when determining damper presence at geometry $g_j$.

An important consideration is the choice of the penalization parameter $\gamma$ during optimization. We address this by employing L-curve plots \cite{Hans01}, which are log-log plots illustrating the trade-off between two quantities that must be balanced. Such curves are common in applied mathematics and engineering, and similar plots have appeared in various studies, see, e.g., \cite{HanOl93, Miller70}. The optimal $\gamma$ is typically found at the \emph{corner} of the L-curve. While L-curves are applicable in various contexts, they are most notably used in Tikhonov regularization (see \cite{Tikhonov1963}), where one minimizes $\|\bA\bx-\bb\|_2^2 + \lambda^2\|\bL(\bx)\|_2^2$ over $\bx$, with $\lambda \in \mathbb{R}$ serving as the regularization parameter chosen by the user.

\subsection{Relaxation via re-weighted $l_1$ minimization}

We begin by considering a relaxed version of the original problem \eqref{eq:optimization_l0}:
\begin{equation}\label{eq:optimization_l1}
	\bar f_\Xi(\vvec)= \trace{\bX(\vvec)} +  \gamma \|\vvec\|_{l_1},
\end{equation}
which can be solved efficiently using convex optimization methods. While $l_1$ relaxation is widely used to approximate $l_0$-norm minimization and has demonstrated considerable potential in various applications (see, e.g., \cite[Sec. 6.5.4 and 11.4.1]{boyd2004convex}), it does face limitations when recovering highly sparse solutions.
To address these limitations, we employ a \emph{weighted $l_1$ relaxation}, as initially proposed in~\cite{CWB08}. Specifically, we consider the following optimization problem:
\begin{equation}\label{eq:optimization_l1_weighted}
	\hf_\Xi(\vvec)= \trace{\bX(\vvec)} +  \gamma \|\vvec\|_{l_1;\wvec},
\end{equation}
where $\wvec = (w_1, \ldots, w_s) \in \mathbb{R}^s_+$ is a vector of nonnegative weights assigned to each component, and the weighted $l_1$-norm is defined as
\begin{equation}
	\|\vvec\|_{l_1;\wvec} = \sum_{j=1}^s w_j \vvec_j.
\end{equation}

The intuition behind the weighted $l_1$ approach is to assign higher weights to smaller values of $\vvec_j$, thereby encouraging them to decrease further, which enhances sparsity by promoting the merging of negligible components. Conversely, larger values of $\vvec_j$ are assigned relatively lower weights, so their contributions remain significant. The main challenge lies in choosing the weighting vector $\wvec$: ideally, it would be the element-wise inverse of $\vvec$, but since $\vvec$ is itself the optimization variable, this is not possible directly.

To circumvent this, we iteratively solve the problem using a re-weighted scheme, as outlined in Algorithm~\ref{alg:solving_opt_iteratively}. At each iteration, the weights are updated based on the current solution, assigning each component a weight inversely proportional to the previous estimate of its viscosity parameter. This iterative procedure was introduced in \cite{CWB08}.

\begin{algorithm}[H]
	\caption{Determining $\vvec_{\mathrm{opt}}$ via iterative re-weighted $l_1$ minimization}
	\label{alg:solving_opt_iteratively}
	\begin{algorithmic}[1]
		\Require The tuple $\Xi$ for $\hf_\Xi(\vvec)$ (e.g., see \eqref{eq:optimization_l1_weighted}), the penalization factor $\gamma > 0$, the number of iterations ($\texttt{iters}$), the initial viscosity vector $\vvec^{(0)} \in \mathbb R^s_+$, initial weights $\wvec^{(0)}$. 
		\Ensure Optimal damping vector $\vvec_{\mathrm{opt}}$.
        \State Set $\varepsilon \gets 10^{-8}$
		\For {$l = 1$ \textbf{to} $\texttt{iters}$}
		\State Use weights $\wvec^{(l-1)}$ to obtain $\vvec^{(l)}$ by minimizing function $\hf_\Xi(\vvec)$ defined in \eqref{eq:optimization_l1_weighted}
		\State Update weights: \quad $w_j^{(l)}=\frac{1}{\vvec_j^{(l)}+\varepsilon}$ for $j=1,\ldots,s$ \label{alg_step:epsilon}
		\EndFor
		\State Set ${\vvec}_{\mathrm{opt}} \gets \vvec^{(l)}$
	\end{algorithmic}	
\end{algorithm}
In Step \ref{alg_step:epsilon} of Algorithm~\ref{alg:solving_opt_iteratively}, the introduction of a small tolerance $\varepsilon > 0$ ensures numerical stability by preventing division by zero. Moreover, it assigns large weights to viscosities that have already attained zero, thereby promoting their persistence at zero in subsequent iterations.

\subsection{Pruning assistance to improve sparsity and estimation of optimal number of dampers}\label{sec:drop_off}

Once we complete the subroutine summarized in \Cref{alg:solving_opt_iteratively}, we obtain $\vvec_{\mathrm{opt}}$, which is the first part of our overall algorithm, as presented later in the paper.

Obtained $\vvec_{\mathrm{opt}}$ may already have zero elements or elements with magnitude up to small tolerance $\texttt{tol}_0$, which means for $|(\vvec_{\mathrm{opt}})_j| < \texttt{tol}_0$ we do not need a damper on position $j$. On the other hand, it is still possible that the number of dampers can be reduced further without a significant increase in function value by dropping \emph{unimportant} dampers. Now the question is how to determine which dampers to drop.

The idea in our approach is to calculate the relative error (relative increase) 
of the function $f_\Xi$ when the damper on each possible position is dropped, one at a time, i.e., for $j \in \{1,\dots, s\}, |(\vvec_{\mathrm{opt}})_j|\geq \texttt{tol}_0$, we calculate
\begin{equation}\label{eq:rel_err_drop}
	\eta_j = \left|\frac{f_\Xi(\vvec_{\mathrm{opt}})-f_\Xi(\vvec_{\mathrm{opt}}^j)}{f_\Xi(\vvec_{\mathrm{opt}})}\right|
\end{equation}
where $(\vvec_{\mathrm{opt}}^j)_j=0$ and $(\vvec_{\mathrm{opt}}^j)_k=(\vvec_{\mathrm{opt}})_k$ for $k \neq j$, i.e., $\vvec_{\mathrm{opt}}^j$ agrees with $\vvec_{\mathrm{opt}}$ except that its $j$th entry is set to zero. Then, the least ``important'' damper is the one whose corresponding relative error $\eta_j$ is the smallest, i.e., we drop the damper on the position $l$ if $\eta_l = \min_{j \in \{1,\dots, s\}}\eta_j $. We repeat this process to determine the ``optimal'' number of dampers denoted by $s^*$ until even the smallest $\eta_j$ is greater than a given tolerance \texttt{tol}. 
This subroutine is summarized in \cref{alg:drop_off_updated}.


\begin{algorithm}[H]
    \caption{Determining optimal number of dampers}
    \label{alg:drop_off_updated}
    \begin{algorithmic}[1]
        \Require Tuple $\Xi$ for $f_\Xi(\vvec)$, optimal viscosity 
        parameter $\vvec_\mathrm{opt} \in \mathbb R^s_+$, and the following scalars: 
        \newline $\texttt{tol}$: relative tolerance for damper elimination,
        \newline $\texttt{tol}_0$: threshold below which a viscosity is deemed inactive.
        \Ensure Optimal number of damping positions $s^*$, and optimal viscosity parameter $\tilde{\vvec}_{\mathrm{opt}}$.
        \State Initialize $\eta^* \gets 0$ and $\tilde{\vvec}_{\mathrm{opt}} \gets \vvec_{\mathrm{opt}}$
        \While{$\eta^* < \texttt{tol}$}
            \State Identify the active index set $\mathcal{A} \gets \{j : |(\tilde{\vvec}_{\mathrm{opt}})_j| \geq \texttt{tol}_0\}$
            \For{each $j \in \mathcal{A}$}
                \State Construct a candidate vector $\tilde{\vvec}_{\mathrm{opt}}^j \gets \tilde{\vvec}_{\mathrm{opt}}$ with $(\tilde{\vvec}_{\mathrm{opt}}^j)_j \gets 0$
                \State Compute the relative change in objective upon removal of damper $j$:
                $$\eta_j \gets \left|\frac{f_\Xi(\vvec_{\mathrm{opt}}) - 
                f_\Xi(\tilde{\vvec}_{\mathrm{opt}}^j)}{f_\Xi(\vvec_{\mathrm{opt}})}\right|$$
            \EndFor
            \State $\eta^* \gets \min_{j \in \mathcal{A}}\, \eta_j$ \Comment{Smallest relative impact} \label{alg_step:etastar}
            \If{$\eta^* < \texttt{tol}$} \Comment{only then is the damper actually removed}
                \State $l \gets \operatorname{argmin}_{j \in \mathcal{A}}\, \eta_j$ \Comment{Index of least influential damper}
                \State Eliminate damper $l$ by setting $(\tilde{\vvec}_{\mathrm{opt}})_l \gets 0$
            \EndIf
        \EndWhile\\
        \Return $s^* \gets |\{j : |(\tilde{\vvec}_{\mathrm{opt}})_j| \geq \texttt{tol}_0\}|$, $\tilde{\vvec}_{\mathrm{opt}}$
    \end{algorithmic}
\end{algorithm}
\begin{remark}\label{rem:adding_dampers} It should be noticed that this approach is different from the approach in which dampers are added consecutively in the system. In the approach where a new damping part is added, it can happen that a good damping position for a smaller number of dampers is no longer good in the case with more damping positions. 

\end{remark}

Finally, we combine all the necessary ingredients needed for the main method, which we call \emph{sparsity-promoting damping optimization} (SPARDO) and which is given in \cref{alg:doil1}. First, in Step \ref{alg:alg_SPARDO_alg1}, we use re-weighted $l_1$ minimization to obtain $\vvec_{\mathrm{opt}}$ as sparse as possible. Second, in Step \ref{st:drop_off}, we prune $\vvec_{\mathrm{opt}}$ with a certain tolerance to obtain $\tilde\vvec_{\mathrm{opt}}$ that has many more zero elements than non-zero elements. In Step \ref{st:deflation} 
we extract positions that contain non-zero elements so we can reduce the dimension of $\tilde\vvec_{\mathrm{opt}}$. Denote the resulting number of non-zero elements by $s^*\in\mathbb{N}$, where $s^*\ll s$. Lastly, we do one more optimization, but with a reduced variable dimension, i.e., we deflate our viscosity vector from dimension $s$ to dimension $s^*$ in Step \ref{st:deflation} and optimize only $s^*$ variables in Step \ref{st:last_min}.  As can be seen, Step \ref{alg:alg_SPARDO_alg1} requires the optimization of $s$ variables. Consequently, the proposed algorithm is computationally feasible only for problems with a moderate dimension of the variable $s$.

\begin{algorithm}[tb]
	\caption{Sparsity-promoting damping optimization (SPARDO)}
	\label{alg:doil1}
	\begin{algorithmic}[1]
		\Require The tuple $\Xi$, the initial viscosity vector $\vvec^{(0)} \in \mathbb R^s_+$, and the following scalars: penalization factor $\gamma > 0$, number of iterations $\texttt{iters}$,  tolerance for determining the number of dampers $\texttt{tol}>0$ , tolerance for whether something is considered zero or not $\texttt{tol}_0$.
		\Ensure Optimal number of dampers in the system $s^*$, their positions $\pi_1,\ldots,\pi_{s^*}$ and viscosities $\tilde\vvec_{\mathrm{opt}}$.
		\State $w_j \gets 1$ (set initial weights for \cref{alg:solving_opt_iteratively} )
        \While {number of dampers did not converge}
		\State $\vvec_{\mathrm{opt}} \gets$ solution returned by \cref{alg:solving_opt_iteratively} using $\Xi$, $\gamma$, $\texttt{iters}$, $\vvec^{(0)}$ and the initial weights. \label{alg:alg_SPARDO_alg1}
		\State $\tilde\vvec_{\mathrm{opt}} \gets$ solution returned by \cref{alg:drop_off_updated} using $\Xi, \vvec_{\mathrm{opt}}$, $\texttt{tol}$ and $\texttt{tol}_0$.  \label{st:drop_off} 
        \State Deflate $\tilde\vvec_{\mathrm{opt}}$ and save positions of its nonzero elements in $\pi_1,\dots,\pi_{s^*}$.\label{st:deflation}
        \EndWhile
        \State $\tilde\vvec_{\mathrm{opt}} \gets \operatorname{argmin}_{\vvec \in\mathbb R^{s^*}_+} \trace{\bX(\vvec)}$ with starting point $\tilde\vvec_{\mathrm{opt}}$, where $\bX$ is the solution of Lyapunov equation \eqref{eq:total_avg_en} where $C_\extern(\vvec)=\sum\limits_{j=1}^{s^*} \vvec_je_{\pi_j}e_{\pi_j}^\T$. \label{st:last_min}
	\end{algorithmic}
	\algcomment{Notice that in Step \ref{st:deflation} the dimension is reduced from $s$ to $s^*$, and in Step \ref{st:last_min} we do one more optimization to refine $\tilde\vvec_{\mathrm{opt}}$.
	}	
\end{algorithm}


\section{Numerical experiments} \label{sec:nume_ex}
\begin{figure}[ht]
  \begin{center}\centering
   \resizebox{.9\linewidth}{!}{
%
%
%

 \begin{tikzpicture}[scale=1, every node/.style={scale=0.8}]
\tikzstyle{spring}=[thick,decorate,decoration={zigzag,pre
length=0.3cm,post length=0.3cm,segment length=6}]
\tikzstyle{damper}=[thick,decoration={markings,
  mark connection node=dmp,
  mark=at position 0.5 with
  {
    \node (dmp) [thick,inner sep=0pt,transform
    shape,rotate=-90,minimum width=10pt,minimum height=2pt,draw=none]
    {};
    \draw [thick] ($(dmp.north east)+(1.5pt,0)$) -- (dmp.south east)
    -- (dmp.south west) -- ($(dmp.north west)+(1.5pt,0)$);
    \draw [thick] ($(dmp.north)+(0,-3pt)$) -- ($(dmp.north)+(0,3pt)$);
  }
}, decorate]
\tikzstyle{springdot}=[thick,decoration={markings,
  mark connection node=sdt,
  mark=at position 0.5 with
  {
  \node (sdt) [thick,inner sep=0pt,transform shape,rotate=-90,minimum
  width=0.85cm,minimum height=0.50cm,draw=none] {};
      \draw [dotted, thick, color=blue] ($(sdt.north)$) --
      ($(sdt.south)$);
  }
}, decorate]

\tikzstyle{ground}=[fill,pattern=north east lines,draw=none,minimum
width=0.75cm,minimum height=0.3cm]

\newcommand{\drawLinewithBG}[2]
{
    \draw[white,myBG]  (#1) -- (#2);
    \draw[black,very thick] (#1) -- (#2);
}

\tikzstyle myBG=[line width=3pt,opacity=1.0]

\node (M) [draw,outer sep=0pt,thick,minimum width=1.2cm, minimum
height=1.2cm,color=red, fill=red!20!white] at (0,0) {$\color{black}
m_1$};
\node (M2) [draw,outer sep=0pt,thick,minimum width=1.2cm, minimum
height=1.2cm,color=red,fill=red!20!white] at (2,0) {$\color{black}
m_2$};
\node (M3) [draw,outer sep=0pt,thick,minimum width=1.2cm, minimum
height=1.2cm,color=red,fill=red!20!white] at (4,0) {$\color{black}
m_{n-1}$};
\node (M4) [draw,outer sep=0pt,thick,minimum width=1.2cm, minimum
height=1.2cm,color=red,fill=red!20!white] at (6,0) {$\color{black}
m_n$};

\node (LW)
[ground,rotate=-90,anchor=north,xshift=-1cm,yshift=-2.0cm,minimum
width=3cm, style={draw,outer sep=0pt,thick}] at (M.south) {};

\node (RW)
[ground,rotate=90,anchor=north,xshift=1cm,yshift=-2.5cm,minimum
width=3cm, style={draw,outer sep=0pt,thick}] at (M4.south west) {};

\draw [spring,color=blue] ({-1.60,0}) -- ({-0.48,0});
\draw [spring,color=blue] (M2.180) -- ($(M.north
east)!(M.180)!(M.south east)$);

\draw [springdot,color=blue] (M3.180) -- ($(M2.north
east)!(M2.180)!(M2.south east)$);
\draw [spring,color=blue] (M4.180) -- ($(M3.north
east)!(M3.180)!(M3.south east)$);
\draw [spring,color=blue] ({6.48,0}) -- ({7.52,0});
\draw [damper] ({0,1.1}) -- ({-1.1,1.1});
\draw [thick] ({0,1.1}) -- ({0,.5});

\draw [damper] ({2,1.1}) -- ({0.9,1.1});
\draw [thick] ({2,1.1}) -- ({2,.5});

\draw [damper] ({4,1.1}) -- ({2.9,1.1});
\draw [thick] ({4,1.1}) -- ({4,.5});

\draw [damper] ({6,1.1}) -- ({4.9,1.1});
\draw [thick] ({6,1.1}) -- ({6,.5});
\node at (-1,-.3) {$k_1$};
\node at (-.7,.7) {$v_1$};

\node at (1,-.3) {$k_2$};

\node at (1.3,.7) {$v_2$};
\node at (3.3,.7) {$v_{n-1}$};
\node at (5.3,.7) {$v_n$};
\node at (5,-.3) {$k_n$};
\node at (7.15,-.3) {$k_{n+1}$};

\node (LW)
[ground,rotate=-90,anchor=north,xshift=-1cm,yshift=-2.0cm,minimum
width=1cm, style={draw,outer sep=0pt,thick}] at (.5,.3) {};
\node (LW)
[ground,rotate=-90,anchor=north,xshift=-1cm,yshift=-2.0cm,minimum
width=1cm, style={draw,outer sep=0pt,thick}] at (2.5,.3) {};
\node (LW)
[ground,rotate=-90,anchor=north,xshift=-1cm,yshift=-2.0cm,minimum
width=1cm, style={draw,outer sep=0pt,thick}] at (4.5,.3) {};
\node (LW)
[ground,rotate=-90,anchor=north,xshift=-1cm,yshift=-2.0cm,minimum
width=1cm, style={draw,outer sep=0pt,thick}] at (6.5,.3) {};

\end{tikzpicture}

%
%
}
  \end{center}
  \caption{$n$-mass oscillator}\label{fig:n_mass}
 \end{figure}
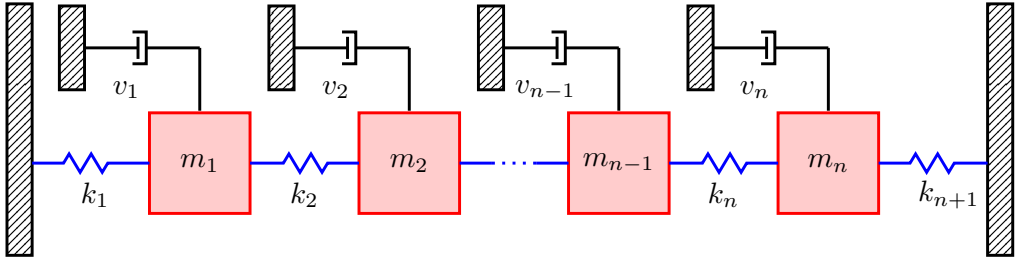

All experiments were performed in MATLAB R2021a using a computer with 2 Intel Xeon Silver 4110 CPUs running at 2.1 GHz and equipped with 192 GB total main memory. Our code for replicating all experiments reported here is both included as supplementary material with this article and available in a permanent and public archive on Zenodo \footnote{Downloadable from https://doi.org/10.5281/zenodo.22204898. We reserved the DOI and the code will be published once the paper is accepted for publication}. As test problems, we use a grounded instance of an $n$-mass oscillator (masses are grounded by a damper); see Figure \ref{fig:n_mass}. One can notice that there is a damper on each mass, which means that the number of dampers $s$ is equal to the dimension of the system $n$. After applying our approach, we will have $s^* \ll n$ which means that we will remove some unnecessary dampers.
	 For this mechanical system, we have the following matrices
\begin{subequations}
\label{eq:nmass}
\setlength\arraycolsep{4pt}
\begin{align}
	M &= \diag{m_1,m_2,\dots,m_n}, \\
    	K &=
	\begin{bsmallmatrix}
	k_1+k_2	& -k_2  	& 			& 			\\
        	-k_2   	& \ddots 	& \ddots   		& 			\\
          		&  \ddots 	& \ddots  		& -k_{n}  		\\
               		&      		&  -k_n 		& k_n+k_{n+1}
        	\end{bsmallmatrix}. 	
\end{align}
\end{subequations}
Here we assume that the number of dampers is initially equal to the dimension of the system, i.e., we optimize over $\vvec=(v_1,\ldots,v_n)\in \mathbb{R}^n_+$, and thus,
\begin{equation}\label{eq:C_ext_n}
C_\extern(\vvec)= \vvec_1 e_{1} e_{1}^\T+\vvec_2 e_{2} e_{2}^\T+\cdots+\vvec_n e_{n} e_{n}^\T,
\end{equation}
where $e_{j}$ for $j=1,\ldots, n$ determine damping geometries, i.e.,  there are grounded dampers on each mass.  We use the Nelder--Mead algorithm \cite{Nelder1965Simplex} implemented in \matlab's routine \texttt{fminsearch} \cite{Lagarias1998NelderMead}. Moreover, for comparison, we use the \texttt{DIRECT} method \cite{Jones1993} from \matlab's toolbox \texttt{DIRECTGO} \cite{DIRECTGO2022} The former is a local method, used with \texttt{TolFun} $=10^{-3}$, \texttt{TolX} $=10^{-2}$ and a budget of $5\,000$ function evaluations for $n=20$ ($50\,000$ for $n=100$), and started from the previous iterate; the latter is a global method, applied on the box $[10^{-13}, 10^{3}]^{s}$ with convergence tolerance $\texttt{tol}_\texttt{DIRECT}$, and is considerably more expensive. Note that the box constraint is imposed only by \texttt{DIRECT}. 
 
 \subsection{Small example}
 In the first example, we illustrate the efficiency of our approach on a small example where we use instances of $n$-mass oscillator defined by the matrices \eqref{eq:nmass}. In this example we use $n = 20$ and define $M$ and $K$ via respectively setting  $m_i = m_{n+1-i} = 10^{-\frac{22}{9}+\frac{4}{9}i}$ for $i=1,\dots,\frac{n}{2},$ and $k_i = 10$ for $i=1,\dots,n+1.$ For damping matrix $C(\vvec)=C_\intern+C_\extern(\vvec)$ we use $\alpha = 0.002$ to define $C_\intern$ given by \eqref{eq:Cint}. Matrix $Z$ needed for the criterion of total average energy is
 \begin{eqnarray*}
  Z &=& \left[\begin{array}{cccc}
                0&0&0&0  \\
                0&Z_1&0&0  \\
                0&0&0&0  \\
                0&0&0&Z_1  \\
  \end{array}\right],
  \end{eqnarray*}
  where 
  $$
  Z(18:20,18:20) = Z(38:40,38:40) = Z_1 = \left[\begin{array}{ccc}
                1&0&0  \\
                0&1&0  \\
                0&0&1 \\
  \end{array}\right], 
  $$
where $Z(18:20,18:20)$ denotes the diagonal block from $(18,18)$ to $(20,20)$. The frequencies we damp with this matrix $Z$ are the following:
\begin{eqnarray*}
    \omega_{18} = 0.5204,\quad \omega_{19} = 0.3315,\quad \omega_{20} = 0.0839.
\end{eqnarray*} 
We use $\vvec = [1,\dots,1]^\T\in \mathbb{R}^{20}_+$ as initial viscosity.  The number of iterations \texttt{iters} in Algorithm \ref{alg:solving_opt_iteratively} for the \texttt{fminsearch} is 20, and for the \texttt{DIRECT} method is 10.
 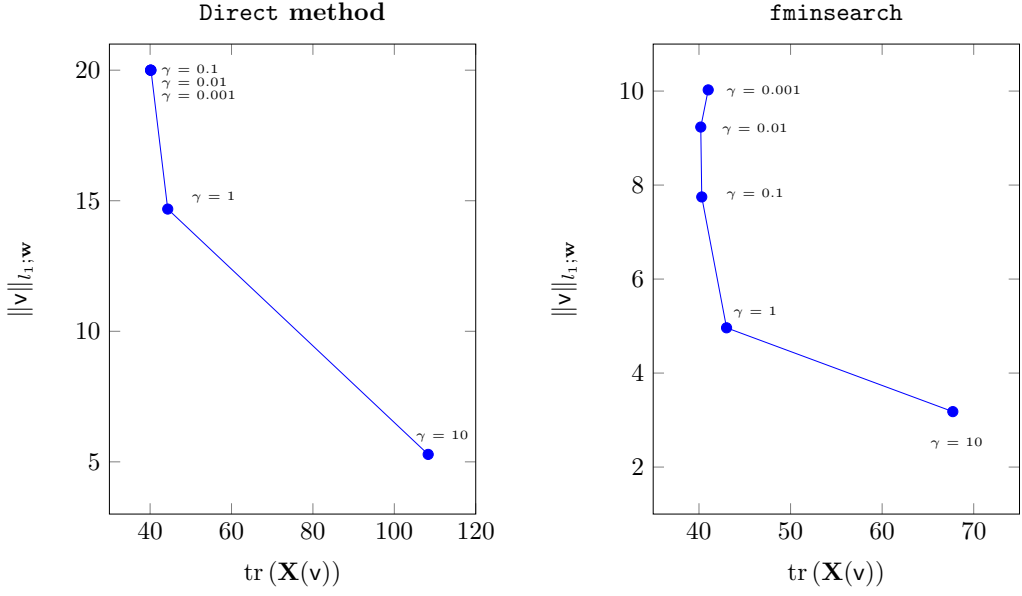
\begin{figure}[t]
  \centering
   \resizebox{0.9\linewidth}{!}{
%
%
\begin{tikzpicture}

\begin{axis}[%
width=2.in,
height=2.566in,
at={(3.527in,0.481in)},
scale only axis,
ymin=1,
ymax=11,
yminorticks=true,
xmin=35,
xmax=75,
xminorticks=true,
axis background/.style={fill=white},
title style={font=\bfseries},
title={\texttt{fminsearch}},
ylabel={$\|\vvec\|_{l_1;\bw}$},
xlabel={$\trace{\bX(\vvec)}$}
]
\addplot [color=blue, mark=*, mark options={solid, blue}]
  table[row sep=crcr]{%
  6.7718e+01   3.1796e+00\\
   4.3011e+01   4.9630e+00\\
   4.0303e+01   7.7474e+00\\
   4.0195e+01   9.2345e+00\\
   4.1007e+01   1.0025e+01\\
};

\node[right, align=left]
at (axis cs:64.26,2.5) {\tiny{$\gamma$ = 10}};
\node[right, align=left]
at (axis cs:42.757,5.3) {\tiny{$\gamma$ = 1}};
\node[right, align=left]
at (axis cs:42,7.8) {\tiny{$\gamma$ = 0.1}};
\node[right, align=left]
at (axis cs:41.5,9.2) {\tiny{$\gamma$ = 0.01}};
\node[right, align=left]
at (axis cs:42, 10) {\tiny{$\gamma$ = 0.001}};
\end{axis}

\begin{axis}[%
width=2in,
height=2.566in,
at={(0.558in,0.481in)},
scale only axis,
ymin=3,
ymax=21,
yminorticks=true,
xmin=30,
xmax=120,
xminorticks=true,
axis background/.style={fill=white},
title style={font=\bfseries},
title={\texttt{Direct} method},
ylabel={$\|\vvec\|_{l_1;\bw}$},
xlabel={$\trace{\bX(\vvec)}$}
]
\addplot [color=blue, mark=*, mark options={solid, blue}]
  table[row sep=crcr]{%
 108.3184  5.2851\\
44.3288 14.6811\\
 40.2149 20\\
40.1240 20\\	
40.1127 20\\
};

\node[right, align=left]
at (axis cs:103.092, 6) {\tiny{$\gamma$ = 10}};
\node[right, align=left]
at (axis cs:48, 15.143) {\tiny{$\gamma$ = 1}};
\node[right, align=left]
at (axis cs:40.5, 20) {\tiny{$\gamma$ = 0.1}};
\node[right, align=left]
at (axis cs:40.5,19.5 ) {\tiny{$\gamma$ = 0.01}};
\node[right, align=left]
at (axis cs:40.5,19) {\tiny{$\gamma$ = 0.001}};
\end{axis}
\end{tikzpicture}%
}
   \caption{The L-curves for our approach with dimension $20$.}\label{fig:l_curve}
 \end{figure}
 
 \cref{fig:l_curve} shows two L-curve plots, one obtained with the \texttt{DIRECT} method and the other with \texttt{fminsearch}. Each curve is traced out by solving \eqref{eq:optimization_l1_weighted} for the penalization parameters $\gamma \in \{10, 1, 10^{-1}, 10^{-2}, 10^{-3}\}$ and plotting, on a log-log scale, the attained weighted $l_1$ norm $\|\vvec_{\mathrm{opt}}\|_{l_1;\wvec}$ against the objective value $\trace{\bX(\vvec_{\mathrm{opt}})}$; the corner of the curve marks the trade-off beyond which a further reduction in the number of dampers costs a disproportionate increase in energy.
 From \cref{fig:l_curve} we see that we should choose between $1$ and $0.1$ for the parameter $\gamma$. Since $\gamma =1$ would result in few dampers in case of \texttt{fminsearch}, we chose a more conservative approach and use $\gamma = 0.1$. 

\begin{figure}[H]
  \centering
   \vspace{-0.3cm}
   \resizebox{0.9\linewidth}{!}{
%
%
\begin{tikzpicture}

\begin{axis}[%
width=4.521in,
height=3.066in,
at={(0.758in,0.481in)},
scale only axis,
xmin=1,
xmax=20,
xlabel style={font=\color{white!15!black}},
xlabel={number of dropped dampers},
ymode=log,
ymin=1e-6,
ymax=10000,
yminorticks=true,
ylabel={relative error $\eta^*$},
ylabel style={font=\color{white!15!black}},
axis background/.style={fill=white},
legend style={legend cell align=left, align=left, draw=white!15!black}
]
\addplot [color=orange]
  table[row sep=crcr]{%
1	0.5\\
2	0.5\\
3	0.5\\
4	0.5\\
5	0.5\\
6	0.5\\
7	0.5\\
8	0.5\\
9	0.5\\
10	0.5\\
11	0.5\\
12	0.5\\
13	0.5\\
14	0.5\\
15	0.5\\
16	0.5\\
17	0.5\\
18	0.5\\
19	0.5\\
20	0.5\\
};
\addlegendentry{0.5}

\addplot [dotted, color=red]
  table[row sep=crcr]{%
1	0.1\\
2	0.1\\
3	0.1\\
4	0.1\\
5	0.1\\
6	0.1\\
7	0.1\\
8	0.1\\
9	0.1\\
10	0.1\\
11	0.1\\
12	0.1\\
13	0.1\\
14	0.1\\
15	0.1\\
16	0.1\\
17	0.1\\
18	0.1\\
19	0.1\\
20	0.1\\
};
\addlegendentry{0.1}

\addplot [dashed, color=green]
  table[row sep=crcr]{%
1	0.01\\
2	0.01\\
3	0.01\\
4	0.01\\
5	0.01\\
6	0.01\\
7	0.01\\
8	0.01\\
9	0.01\\
10	0.01\\
11	0.01\\
12	0.01\\
13	0.01\\
14	0.01\\
15	0.01\\
16	0.01\\
17	0.01\\
18	0.01\\
19	0.01\\
20	0.01\\
};
\addlegendentry{0.01}

\addplot [color=blue, draw=none, only marks, mark=*, mark options={solid, blue}]
  table[row sep=crcr]{%
1     1.668983003482271e-06 \\
2     3.346952332039030e-06\\
3     1.011266056266341e-05\\
4     1.691418290202246e-05\\
5     3.263576502896859e-05\\
6     4.843754987911096e-05\\
7     7.783000141268296e-05\\
8     1.073642901334019e-04\\
9     1.478585547591125e-04\\
10     1.885727195749218e-04\\
11    2.422894765463219e-04\\
12     2.963176085053327e-04\\
13     9.725701528404688e-03\\
14     1.938007630749244e-02\\
15     9.592135734235550e-02\\
16     2.274375460349812e-01\\
17     4.473476689465512e-01\\
18     1.179625521882414e+00\\
19     3.721783724807314e+00\\
20     4.180322779140997e+02\\
};

\end{axis}

\begin{axis}[%
width=5.833in,
height=3.375in,
at={(0in,0in)},
scale only axis,
xmin=0,
xmax=1,
ymin=0,
ymax=1,
axis line style={draw=none},
ticks=none,
axis x line*=bottom,
axis y line*=left,
legend style={legend cell align=left, align=left, draw=white!15!black}
]
\end{axis}
\end{tikzpicture}%
}
      \caption{Behavior of relative error $\eta^*$ (computed in \Cref{alg:drop_off_updated}, Step~\ref{alg_step:etastar}) in function value as we drop one damper at a time, where the optimal viscosity $\vvec_{\mathrm{opt}}$ is computed with \texttt{DIRECT} method.} \label{fig:damper_drop_off_direct}
 \end{figure}
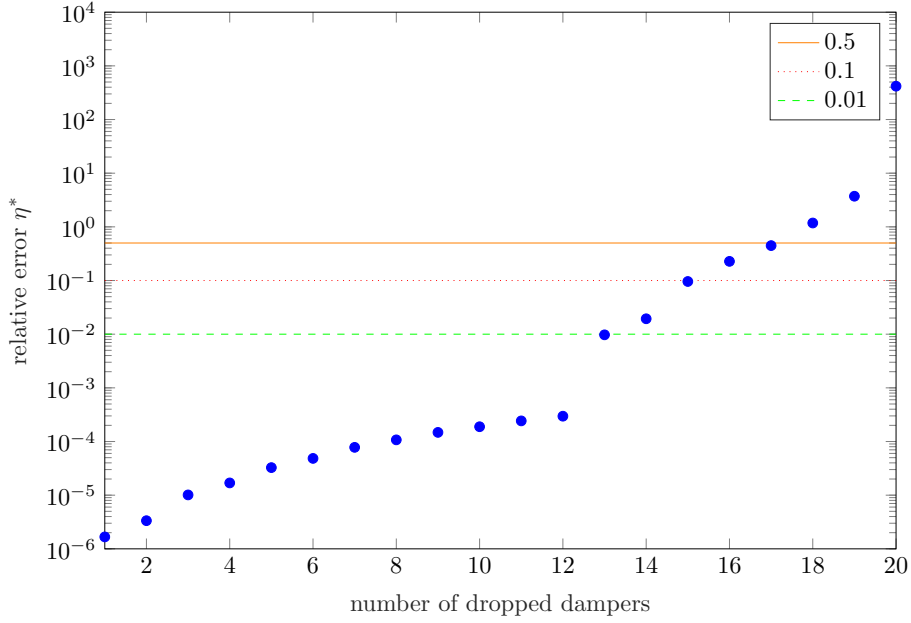

\begin{figure}[H]
  \centering
 \vspace{-0.3cm}
   \resizebox{0.9\linewidth}{!}{
%
%
\begin{tikzpicture}

\begin{axis}[%
width=4.521in,
height=3.066in,
at={(0.758in,0.481in)},
scale only axis,
xmin=1,
xmax=20,
xlabel style={font=\color{white!15!black}},
xlabel={number of dropped dampers},
ymode=log,
ymin=1e-16,
ymax=10000,
yminorticks=true,
ylabel={relative error $\eta^*$},
ylabel style={font=\color{white!15!black}},
axis background/.style={fill=white},
legend style={legend cell align=left, align=left, draw=white!15!black}
]
\addplot [color=orange]
  table[row sep=crcr]{%
1	0.5\\
2	0.5\\
3	0.5\\
4	0.5\\
5	0.5\\
6	0.5\\
7	0.5\\
8	0.5\\
9	0.5\\
10	0.5\\
11	0.5\\
12	0.5\\
13	0.5\\
14	0.5\\
15	0.5\\
16	0.5\\
17	0.5\\
18	0.5\\
19	0.5\\
20	0.5\\
};
\addlegendentry{0.5}

\addplot [dotted, color=red]
  table[row sep=crcr]{%
1	0.1\\
2	0.1\\
3	0.1\\
4	0.1\\
5	0.1\\
6	0.1\\
7	0.1\\
8	0.1\\
9	0.1\\
10	0.1\\
11	0.1\\
12	0.1\\
13	0.1\\
14	0.1\\
15	0.1\\
16	0.1\\
17	0.1\\
18	0.1\\
19	0.1\\
20	0.1\\
};
\addlegendentry{0.1}

\addplot [dashed, color=green]
  table[row sep=crcr]{%
1	0.01\\
2	0.01\\
3	0.01\\
4	0.01\\
5	0.01\\
6	0.01\\
7	0.01\\
8	0.01\\
9	0.01\\
10	0.01\\
11	0.01\\
12	0.01\\
13	0.01\\
14	0.01\\
15	0.01\\
16	0.01\\
17	0.01\\
18	0.01\\
19	0.01\\
20	0.01\\
};
\addlegendentry{0.01}

\addplot [color=blue, draw=none, only marks, mark=*, mark options={solid, blue}]
  table[row sep=crcr]{%
1     2.450142140682014e-15\\
2     9.975578715633915e-15\\
3    5.250304587175745e-16\\
4     4.900284281364028e-15\\
5     2.450142140682014e-15\\
6     2.642653308878458e-14\\
7     4.235245700321767e-14\\
8     5.142320237446502e-02\\
9     1.162239631844166e-01\\
10     2.953889159775512e-01\\
11     7.803117007642462e-01\\
12     2.898892809986723e+00\\
13    4.140564067355322e+02\\
};
\end{axis}

\begin{axis}[%
width=5.833in,
height=3.375in,
at={(0in,0in)},
scale only axis,
xmin=0,
xmax=1,
ymin=0,
ymax=1,
axis line style={draw=none},
ticks=none,
axis x line*=bottom,
axis y line*=left,
legend style={legend cell align=left, align=left, draw=white!15!black}
]
\end{axis}
\end{tikzpicture}%
}
      \caption{Behavior of relative error $\eta^*$ (computed in \Cref{alg:drop_off_updated}, Step~\ref{alg_step:etastar}) in function value as we drop one damper at a time, where the optimal viscosity $\vvec_{\mathrm{opt}}$ is computed with \texttt{fminsearch}. } \label{fig:damper_drop_off}
 \end{figure}
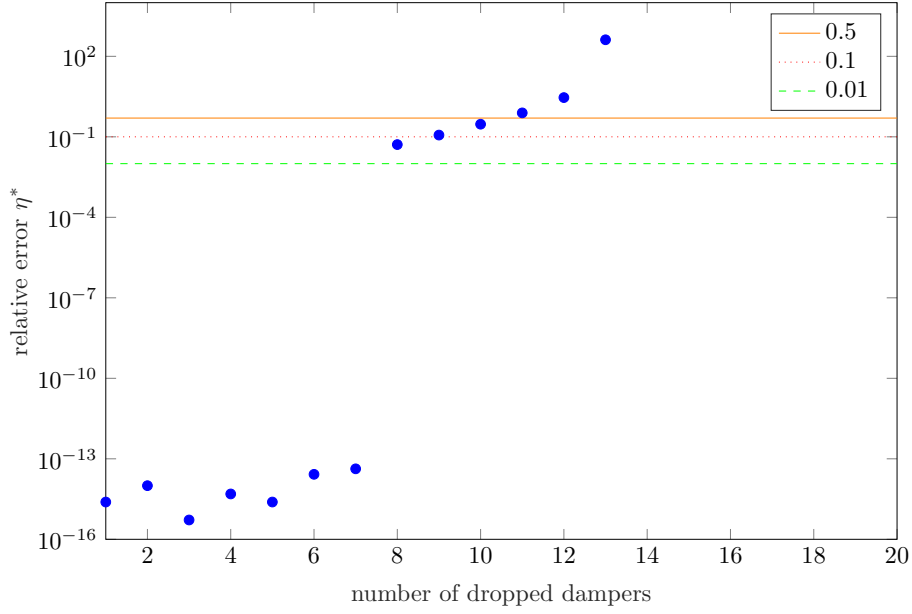

\cref{fig:damper_drop_off_direct} and \cref{fig:damper_drop_off} show the change in relative error $\eta^*$, computed in \Cref{alg:drop_off_updated}, Step~\ref{alg_step:etastar}, with respect to number of dropped dampers. The order in which we drop dampers is explained in \cref{sec:drop_off}. In \cref{fig:damper_drop_off_direct} the optimal viscosity $\vvec_{\mathrm{opt}}$ is computed with \texttt{DIRECT} method. Additionally, from \cref{fig:damper_drop_off_direct}, one can see that in this case $\|\vvec_{\mathrm{opt}}\|_{l_0}=20$, which means we were unable to detect unimportant dampers in Algorithm \ref{alg:doil1} Step \ref{alg:alg_SPARDO_alg1}, that is,  before pruning, which is done in Algorithm \ref{alg:doil1} Step \ref{st:drop_off}. The orange line, red dotted line, and green dashed line represent, respectively, 0.5, 0.1, 0.01 tolerance \texttt{tol} in \cref{alg:drop_off_updated}. This indicates that using at least seven dampers ensures that the relative increase in the objective function value does not exceed 0.01,  i.e., $\eta^*\leq 0.01$. Similarly, at least five dampers are required when $\texttt{tol} =0.1$ and three dampers are enough for $\texttt{tol} =0.5$. 

 On the other hand, in \cref{fig:damper_drop_off} the optimal viscosity $\vvec_{\mathrm{opt}}$ is computed with \texttt{fminsearch} and one can see that in this case $\|\vvec_{\mathrm{opt}}\|_{l_0}=13$, which means that we initially, in Step \ref{alg:alg_SPARDO_alg1} of Algorithm \ref{alg:doil1}, dropped 7 unimportant dampers. The orange line, red dotted line, and green dashed line again represent, respectively, 0.5, 0.1, and 0.01 tolerance \texttt{tol} in \cref{alg:drop_off_updated}.  In this case, using at least six dampers ensures that the relative increase in the objective function value does not exceed 0.01,  i.e., $\eta^*\leq 0.01$, while for other two tolerances we require the same number of dampers as in \cref{fig:damper_drop_off_direct} where the computation is done with \texttt{DIRECT} method, that is, at least five dampers are required when $\texttt{tol} =0.1$ and 3 dampers are enough for $\texttt{tol} =0.5$. 
 
By the ``brute force'' approach we mean the exhaustive search in which, for a prescribed number of dampers, the viscosities are optimized separately for \emph{every} admissible set of positions and the configuration with the smallest function value is retained. If we were to determine both the number of dampers and their positions in this way, we would have to optimize the viscosities for all $2^{20}-1 = 1\,048\,575$ nonempty subsets of the $20$ candidate positions and choose the one with the smallest function value. Even for a fixed number of dampers the counts grow quickly, namely $\binom{20}{3} = 1\,140$, $\binom{20}{5} = 15\,504$, $\binom{20}{6} = 38\,760$ and $\binom{20}{7} = 77\,520$, so that the four cases considered here already amount to $132\,924$ configurations. Thus, for comparison, we only computed optimal viscosities and positions for 3, 5, 6, and 7 dampers.  Obtained results are shown in Tables \ref{tab:Table7dim20} - \ref{tab:Table3dim20tol6} where we show the damper positions, the corresponding objective function values \eqref{eq:min_f}, and the time needed for the computation, obtained for different numbers of dampers using the \texttt{DIRECT} and \texttt{fminsearch} methods. We need to emphasize that the time needed for creating the L-curve is not included in the optimization time.

 \begin{table}[H]
  \begin{center}
   \begin{tabular}{l|c|c|c|c}
   &\multicolumn{2}{c|}{\texttt{DIRECT} method}&\multicolumn{2}{c}{\texttt{fminsearch} method}\\
	& ``brute force'' & SPARDO &``brute force'' & SPARDO \\ \hline
	Optimization time &	$5.7798\cdot 10^5$ & 963.5454 & $1.7362\cdot 10^5$&  97.3226\\
	Optimal position & (7:13)
 & (7:13) &(7:13)
 & [6,(8:13)] \\
	Function value & $40.3812$ & $40.3812$ & $40.3976$ & $40.5951$
   \end{tabular}
  \end{center}
     \vspace{-0.1cm}  
  \caption{Results for 7 dampers,  $\gamma=0.1$, $\texttt{tol}_\texttt{DIRECT}=10^{-3}$  }
  \label{tab:Table7dim20}
 \end{table}
 Table \ref{tab:Table7dim20} presents results for seven dampers. When the \texttt{DIRECT} method was employed, SPARDO (Algorithm \ref{alg:doil1}) produced the same damper positions and the identical objective function value as the brute-force approach, while reducing the computational time by nearly a factor of 600.

In contrast, when the \texttt{fminsearch} method was used, the relative error between the objective function values obtained by the ``brute force'' approach and SPARDO was $4.8889 \cdot 10^{-3}$, although SPARDO was more than 1700 times faster. This discrepancy arises because in this case the SPARDO algorithm identifies six as the optimal number of dampers, as illustrated in \cref{fig:damper_drop_off}. Consequently, for the case of seven dampers, the most accurate result is obtained when SPARDO is combined with the \texttt{DIRECT} method.

  \begin{table}[H]
  \begin{center}
   \begin{tabular}{l|c|c|c|c}
   &\multicolumn{2}{c|}{\texttt{DIRECT} method}&\multicolumn{2}{c}{\texttt{fminsearch} method}\\
	& ``brute force'' & SPARDO &``brute force'' & SPARDO \\ \hline
	Optimization time &	$2.1264\cdot 10^{5}$ & 960.5135&$ 7.0523\cdot 10^{4}$ &  97.3092\\
	Optimal position & (8:13)
 & (8:13) &(7,9,10,11,12,13)  &(8:13) \\
	Function value & $40.6097$ & $40.6097 $ & $40.8205$ &$40.5951 $
   \end{tabular}
  \end{center}
   \vspace{-0.1cm}
  \caption{Results for 6 dampers, $\gamma=0.1$, $\texttt{tol}_\texttt{DIRECT}=10^{-3}$}
     \label{tab:Table6dim20}
 \end{table}
 
 \vspace{-0.5cm}
  \begin{table}[H]
  \begin{center}
   \begin{tabular}{l|c|c|c|c}
   &\multicolumn{2}{c|}{\texttt{DIRECT} method}&\multicolumn{2}{c}{\texttt{fminsearch} method}\\
	& ``brute force'' & SPARDO &``brute force'' & SPARDO \\ \hline
	Optimization time &	$6.7246\cdot 10^{4}$ & 958.9407&$2.2341\cdot 10^{4}$ &  97.4569\\
	Optimal position & (8,9,10,11,13) & (8,9,11,12,13) &(8,10,11,12,13)  &(8,10,11,12,13) \\
	Function value & $42.0867$ & $42.8128 $ & $42.0818$ &$42.0818$ 
   \end{tabular}
  \end{center}
   \vspace{-0.1cm}
  \caption{Results for 5 dampers, $\gamma=0.1$, $\texttt{tol}_\texttt{DIRECT}=10^{-3}$}
     \label{tab:Table5dim20}
 \end{table}

 \cref{tab:Table6dim20} and \cref{tab:Table5dim20} show similar results. However, compared to \cref{tab:Table7dim20}, the best result is obtained when SPARDO is combined with the \texttt{fminsearch} method. 
 
 For six dampers, the relative error in function value between ``brute force'' and the SPARDO algorithm is $5.5795 \cdot 10^{-3}$ when the \texttt{fminsearch} method is used, whereas the \texttt{DIRECT} method yields the same function values. In both cases, SPARDO is more than 200 times faster than the ``brute force'' approach. 
 On the other hand, for five dampers, the relative error in function value between ``brute force'' and SPARDO is $1.6960 \cdot 10^{-2}$ when the \texttt{DIRECT} method is used, while we obtained the same function values with the \texttt{fminsearch} method. Moreover, SPARDO again provides a substantial reduction in computational time for both optimization methods.
 
     \begin{table}[H]
  \begin{center}
   \begin{tabular}{l|c|c|c|c}
   &\multicolumn{2}{c|}{\texttt{DIRECT} method}&\multicolumn{2}{c}{\texttt{fminsearch} method}\\
	& ``brute force'' & SPARDO &``brute force'' & SPARDO \\ \hline
	Optimization time &	$3.2337\cdot 10^{3}$  & 956.5514 & 853.2421 & 97.4001\\
	Optimal position & (8,11,13) & (8,11,13) & (8,11,13) &(8,10,13) \\
	Function value & $47.5983$ & $47.5983$ & $47.5982$ &$47.5982$
   \end{tabular}
  \end{center}
     \vspace{-0.1cm}
  \caption{Results for 3 dampers, $\gamma=0.1$, $\texttt{tol}_\texttt{DIRECT}=10^{-3}$}
  \label{tab:Table3dim20}
 \end{table}

The results for three dampers are presented in \cref{tab:Table3dim20}. It can be seen that all four approaches produce almost identical results. Using the \texttt{DIRECT} method, SPARDO yields the same damper positions and objective function value as the ``brute force'' approach while being more than three times faster. In contrast, when the \texttt{fminsearch} method is used, the damper positions obtained by SPARDO differ slightly from those of the ``brute force'' approach. However, the relative error in the objective function value is of the order of $10^{-13}$, indicating the existence of two local minima with essentially identical objective function values. In this case, SPARDO is more than eight times faster than the ``brute force'' approach.

A difference in the objective function values obtained with the \texttt{DIRECT} and \texttt{fminsearch} methods can also be observed. This difference is eliminated by using a stricter tolerance in the \texttt{DIRECT} method, and the results are shown in \cref{tab:Table3dim20tol6}. The results presented in Tables \ref{tab:Table7dim20}--\ref{tab:Table3dim20} were obtained with the tolerance $\texttt{tol}_{\texttt{DIRECT}} = 10^{-3}$, whereas the results in \cref{tab:Table3dim20tol6} were obtained with $\texttt{tol}_{\texttt{DIRECT}} = 10^{-6}$. As shown in \cref{tab:Table3dim20tol6}, the stricter tolerance produces the same objective function value for all four approaches, up to an error of the order of $10^{-13}$, while SPARDO remains computationally more efficient than the brute-force approach.

    \begin{table}[H]
  \begin{center}
   \begin{tabular}{l|c|c|c|c}
   &\multicolumn{2}{c|}{\texttt{DIRECT} method}&\multicolumn{2}{c}{\texttt{fminsearch} method}\\
	& ``brute force'' & SPARDO &``brute force'' & SPARDO \\ \hline
	Optimization time &	$9.5451\cdot 10^5$  &  $1.9548\cdot 10^4$ &  853.2421 &97.4001\\
	Optimal position & (8,11,13) & (8,10,13) &(8,11,13) &(8,10,13) \\
	Function value & $47.5982$ & $47.5982$ & $47.5982$ &$47.5982$ 
   \end{tabular}
  \end{center}
     \vspace{-0.1cm}
  \caption{Results for 3 dampers $\gamma=0.1$, $\texttt{tol}_\texttt{DIRECT}=10^{-6}$}
  \label{tab:Table3dim20tol6}
 \end{table}

Based on the results presented in the tables, SPARDO combined with the \texttt{fminsearch} method offers the best overall performance, yielding the lowest objective function values while requiring the least computational time.

 \begin{remark}
Regarding the Remark \ref{rem:adding_dampers} under \cref{alg:drop_off_updated}, we also attempted to determine the position of a single damper from which to initiate the approach of consecutive damper addition. The result was either position $9$ or position $12$, neither of which is in the optimal set of 3 dampers.\end{remark}

 \subsection{Moderate example} \label{exam:moderate}
 In this example we show the results for a similar example of $n$-mass oscillator of dimension $n=100$. We defined matrices $M$ and $K$ by setting $m_i$ = $m_{n+1-i} = 10^{-\frac{102}{49}+\frac{4}{49}i}$ for $i = 1, \dots, \frac{n}{2}$  and $k_i=10$ for $i = 1,\dots, n+1$, while the damping matrix is the same as in the previous example, i.e. $C(\vvec)=C_\intern+C_\extern(\vvec)$ with $\alpha = 0.002$ that determines  $C_\intern$ in \eqref{eq:Cint}.  The matrix $Z$ needed for the criterion of total average energy is
 \begin{eqnarray*}
  Z &=& \left[\begin{array}{cccc}
                0&0&0&0  \\
                0&Z_1&0&0  \\
                0&0&0&0  \\
                0&0&0&Z_1  \\
  \end{array}\right],
  \end{eqnarray*}
  where 
  $$
  Z(99:100,99:100) = Z(199:200,199:200) = Z_1 = \left[\begin{array}{cc}
                1&0  \\
                0&1  \\
  \end{array}\right], 
  $$ 
where $Z(99:100,99:100)$ denotes the diagonal block from $(99,99)$ to $(100,100)$. The frequencies we damp with this matrix $Z$ are the following:
\begin{eqnarray*}
    \omega_{99} = 0.0790,\quad \omega_{100} = 0.0200.
\end{eqnarray*}
We use $\vvec = [100,\dots,100]^\T\in \mathbb{R}^{100}_+$ as initial viscosity. 

In this example, we decided to test our approach only with \texttt{fminsearch}, since the \texttt{DIRECT} method was rather slow even for a small dimension. The number of iterations we used in Algorithm \ref{alg:solving_opt_iteratively} for \texttt{fminsearch} is 30.
 \begin{figure}[H]
  \centering
   \resizebox{0.5\linewidth}{!}{
%
%
\begin{tikzpicture}

\begin{axis}[%
width=2.in,
height=2.566in,
at={(3.527in,0.481in)},
scale only axis,
ymin=58,
ymax=85,
yminorticks=true,
xmin=135,
xmax=150,
xminorticks=true,
axis background/.style={fill=white},
title style={font=\bfseries},
title={\texttt{fminsearch}},
ylabel={$\|\vvec\|_{l_1;\bw}$},
xlabel={$\trace{\bX(\vvec)}$}
]
\addplot [color=blue, mark=*, mark options={solid, blue}]
  table[row sep=crcr]{%
147.600044138649	71.2105241899249\\
137.132973358566	70.9965949968629\\
135.777332538471	60.0589578148481\\
135.533156277661	72.6234889193449\\
135.476620121116	84.3483070981961\\
};

\node[right, align=left]
at (axis cs:145.6,72) {\footnotesize{$\gamma$ = 10}};
\node[right, align=left]
at (axis cs:137.133,70) {\footnotesize{$\gamma$ = 1}};
\node[right, align=left]
at (axis cs:135.777,60.059) {\footnotesize{$\gamma$ = 0.1}};
\node[right, align=left]
at (axis cs:135.533,72.623) {\footnotesize{$\gamma$ = 0.01}};
\node[right, align=left]
at (axis cs:135.477,84) {\footnotesize{$\gamma$ = 0.001}};
\end{axis}

\end{tikzpicture}%
}
   \caption{The L curve for our approach with dimension $100$.}\label{fig:l_curve100}
 \end{figure}
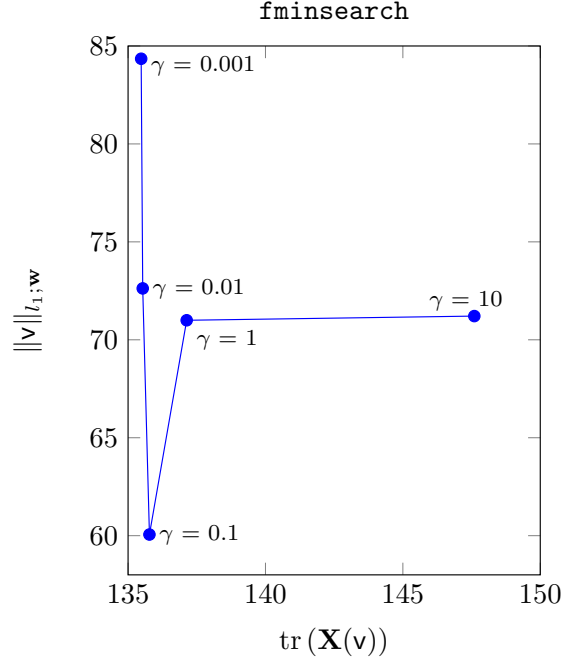
 From \cref{fig:l_curve100}, we see that we should again use $0.1$ for the parameter $\gamma$. Note, however, that although $\|\vvec\|_{l_1;\wvec}=60$ at this value of $\gamma$, the re-weighted $l_1$ stage does not by itself render any of the viscosities inactive, i.e., $\|\vvec\|_{l_0}=100$, so the reduction is left entirely to the pruning stage.

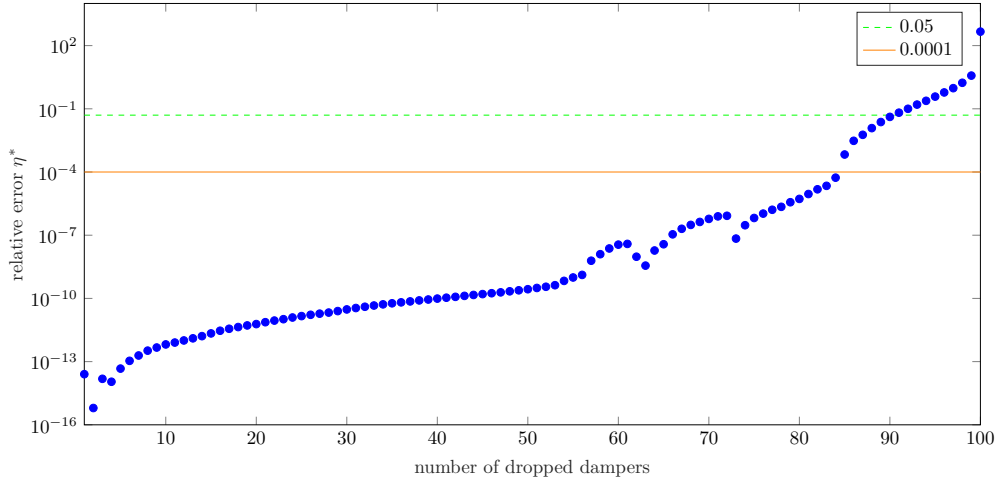
\begin{figure}[H]
  \centering
   \vspace{-0.3cm}
   \resizebox{0.9\linewidth}{!}{
%
%
\begin{tikzpicture}

\begin{axis}[%
width=6.521in,
height=3.066in,
at={(0.758in,0.481in)},
scale only axis,
xmin=1,
xmax=100,
xlabel style={font=\color{white!15!black}},
xlabel={number of dropped dampers},
ymode=log,
ymin=1e-16,
ymax=10000,
yminorticks=true,
ylabel={relative error $\eta^*$},
ylabel style={font=\color{white!15!black}},
axis background/.style={fill=white},
legend style={legend cell align=left, align=left, draw=white!15!black}
]
\addplot [dashed, color=green]
  table[row sep=crcr]{%
1	0.05\\
2	0.05\\
3	0.05\\
4	0.05\\
5	0.05\\
6	0.05\\
7	0.05\\
8	0.05\\
9	0.05\\
10	0.05\\
11	0.05\\
12	0.05\\
13	0.05\\
14	0.05\\
15	0.05\\
16	0.05\\
17	0.05\\
18	0.05\\
19	0.05\\
20	0.05\\
21	0.05\\
22	0.05\\
23	0.05\\
24	0.05\\
25	0.05\\
26	0.05\\
27	0.05\\
28	0.05\\
29	0.05\\
30	0.05\\
31	0.05\\
32	0.05\\
33	0.05\\
34	0.05\\
35	0.05\\
36	0.05\\
37	0.05\\
38	0.05\\
39	0.05\\
40	0.05\\
41	0.05\\
42	0.05\\
43	0.05\\
44	0.05\\
45	0.05\\
46	0.05\\
47	0.05\\
48	0.05\\
49	0.05\\
50	0.05\\
51	0.05\\
52	0.05\\
53	0.05\\
54	0.05\\
55	0.05\\
56	0.05\\
57	0.05\\
58	0.05\\
59	0.05\\
60	0.05\\
61	0.05\\
62	0.05\\
63	0.05\\
64	0.05\\
65	0.05\\
66	0.05\\
67	0.05\\
68	0.05\\
69	0.05\\
70	0.05\\
71	0.05\\
72	0.05\\
73	0.05\\
74	0.05\\
75	0.05\\
76	0.05\\
77	0.05\\
78	0.05\\
79	0.05\\
80	0.05\\
81	0.05\\
82	0.05\\
83	0.05\\
84	0.05\\
85	0.05\\
86	0.05\\
87	0.05\\
88	0.05\\
89	0.05\\
90	0.05\\
91	0.05\\
92	0.05\\
93	0.05\\
94	0.05\\
95	0.05\\
96	0.05\\
97	0.05\\
98	0.05\\
99	0.05\\
100	0.05\\
};
\addlegendentry{0.05}

\addplot [color=orange]
  table[row sep=crcr]{%
1	0.0001\\
2	0.0001\\
3	0.0001\\
4	0.0001\\
5	0.0001\\
6	0.0001\\
7	0.0001\\
8	0.0001\\
9	0.0001\\
10	0.0001\\
11	0.0001\\
12	0.0001\\
13	0.0001\\
14	0.0001\\
15	0.0001\\
16	0.0001\\
17	0.0001\\
18	0.0001\\
19	0.0001\\
20	0.0001\\
21	0.0001\\
22	0.0001\\
23	0.0001\\
24	0.0001\\
25	0.0001\\
26	0.0001\\
27	0.0001\\
28	0.0001\\
29	0.0001\\
30	0.0001\\
31	0.0001\\
32	0.0001\\
33	0.0001\\
34	0.0001\\
35	0.0001\\
36	0.0001\\
37	0.0001\\
38	0.0001\\
39	0.0001\\
40	0.0001\\
41	0.0001\\
42	0.0001\\
43	0.0001\\
44	0.0001\\
45	0.0001\\
46	0.0001\\
47	0.0001\\
48	0.0001\\
49	0.0001\\
50	0.0001\\
51	0.0001\\
52	0.0001\\
53	0.0001\\
54	0.0001\\
55	0.0001\\
56	0.0001\\
57	0.0001\\
58	0.0001\\
59	0.0001\\
60	0.0001\\
61	0.0001\\
62	0.0001\\
63	0.0001\\
64	0.0001\\
65	0.0001\\
66	0.0001\\
67	0.0001\\
68	0.0001\\
69	0.0001\\
70	0.0001\\
71	0.0001\\
72	0.0001\\
73	0.0001\\
74	0.0001\\
75	0.0001\\
76	0.0001\\
77	0.0001\\
78	0.0001\\
79	0.0001\\
80	0.0001\\
81	0.0001\\
82	0.0001\\
83	0.0001\\
84	0.0001\\
85	0.0001\\
86	0.0001\\
87	0.0001\\
88	0.0001\\
89	0.0001\\
90	0.0001\\
91	0.0001\\
92	0.0001\\
93	0.0001\\
94	0.0001\\
95	0.0001\\
96	0.0001\\
97	0.0001\\
98	0.0001\\
99	0.0001\\
100	0.0001\\
};
\addlegendentry{0.0001}

\addplot [color=blue, draw=none, only marks, mark=*, mark options={solid, blue}]
  table[row sep=crcr]{%
1	2.52990137411147e-14\\
2	6.27248274573092e-16\\
3	1.52630413479452e-14\\
4	1.10813861841246e-14\\
5	4.66254550765998e-14\\
6	1.08932117017527e-13\\
7	1.95074213392232e-13\\
8	3.31814337249166e-13\\
9	4.66045468007807e-13\\
10	6.54010867621544e-13\\
11	8.02250543178985e-13\\
12	1.0094515565463e-12\\
13	1.27352308014157e-12\\
14	1.617673300124e-12\\
15	2.18115133344883e-12\\
16	2.91984071813774e-12\\
17	3.62403144772513e-12\\
18	4.36878423240159e-12\\
19	5.21180591342782e-12\\
20	6.04144629792983e-12\\
21	7.42390149508893e-12\\
22	8.87012693349628e-12\\
23	1.0327642840846e-11\\
24	1.24055072917484e-11\\
25	1.44779355909379e-11\\
26	1.66843859381279e-11\\
27	1.88912544508342e-11\\
28	2.12785613838594e-11\\
29	2.50389147899251e-11\\
30	3.00023303866219e-11\\
31	3.49989901418712e-11\\
32	4.01781791450212e-11\\
33	4.5974371367593e-11\\
34	5.22238550099229e-11\\
35	5.86075697830114e-11\\
36	6.50084293422716e-11\\
37	7.20934076863329e-11\\
38	8.03814482210254e-11\\
39	8.86251632109813e-11\\
40	9.80259421847664e-11\\
41	1.08048951446928e-10\\
42	1.18684991356106e-10\\
43	1.31895049101373e-10\\
44	1.45828324107224e-10\\
45	1.6024290762243e-10\\
46	1.76189231424143e-10\\
47	1.93865296884371e-10\\
48	2.17193078381261e-10\\
49	2.42475156418962e-10\\
50	2.73046191907864e-10\\
51	3.13381810369803e-10\\
52	3.5720472015561e-10\\
53	4.23336969323158e-10\\
54	6.79505373741567e-10\\
55	9.84496065776889e-10\\
56	1.30532791298655e-09\\
57	6.15447727333893e-09\\
58	1.25430013680314e-08\\
59	2.34440936080231e-08\\
60	3.56848553949514e-08\\
61	3.87043861449464e-08\\
62	9.52742960006279e-09\\
63	3.6013495228683e-09\\
64	1.87860170352367e-08\\
65	3.72266337446797e-08\\
66	1.09952270265969e-07\\
67	2.01691368855758e-07\\
68	3.08636062678431e-07\\
69	4.26669492362345e-07\\
70	5.94849494273512e-07\\
71	7.84464684389858e-07\\
72	8.3540793928644e-07\\
73	6.90685457441696e-08\\
74	2.97505522810518e-07\\
75	6.5742120603626e-07\\
76	1.06072662598921e-06\\
77	1.61220010668235e-06\\
78	2.22516295965952e-06\\
79	3.67359598436069e-06\\
80	5.28687726195822e-06\\
81	9.0099525127332e-06\\
82	1.50921132224762e-05\\
83	2.2051174123241e-05\\
84	5.38884457675946e-05\\
85	0.000676295953812928\\
86	0.00302653078604042\\
87	0.00582446758374957\\
88	0.0121086467396933\\
89	0.0235355552951202\\
90	0.0413611004500051\\
91	0.0654765952633822\\
92	0.10018358369275\\
93	0.158432952809719\\
94	0.239396774551748\\
95	0.379276268278818\\
96	0.597668814318493\\
97	0.954656534003225\\
98	1.7194130904044\\
99	3.76756043585074\\
100	460.298894898803\\
};
\end{axis}

\begin{axis}[%
width=5.833in,
height=3.375in,
at={(0in,0in)},
scale only axis,
xmin=0,
xmax=1,
ymin=0,
ymax=1,
axis line style={draw=none},
ticks=none,
axis x line*=bottom,
axis y line*=left,
legend style={legend cell align=left, align=left, draw=white!15!black}
]
\end{axis}
\end{tikzpicture}%
}
      \caption{Behavior of relative error $\eta^*$ (computed in \Cref{alg:drop_off_updated}, Step~\ref{alg_step:etastar}) in function value as we drop one damper at a time, where the optimal viscosity $\vvec_{\mathrm{opt}}$ is computed with \texttt{fminsearch} method for $n=100$.} \label{fig:damper_drop_off_100}
 \end{figure}

 \cref{fig:damper_drop_off_100} shows the change in relative error $\eta^*$, compared to the number of dropped dampers. The order in which we drop the dampers is explained in \cref{sec:drop_off}. The orange line and green dashed line represent, respectively, 0.0001, 0.05, tolerance \texttt{tol} in \cref{alg:drop_off_updated}, which means that for at most 16 dampers the relative increase in function value will not exceed 0.0001, i.e., $\eta^*\leq 0.0001$. In the same way, we need 10 dampers for $\texttt{tol} =0.05$.  
 In this example, it is impossible to make a comparison with the ``brute force'' approach, since we would need to check $100\choose 10$ combinations, which would take over $10^5$ years if we assume that we need only 2 seconds for optimizing one combination. Thus, we present only our results. The optimal positions and viscosities of the dampers are summarized in \cref{tab:Example2_pos}. The corresponding value of the objective function is $136.6032$, while the time needed for this computation is $2.5083 \cdot 10^4$ seconds (almost 7 hours).
 \begin{table}[H]
  \begin{center}
   \begin{tabular}{c|c}
	Positions & Viscosities  \\ \hline
   37&    6.7406\\
   41&    4.5082\\
   43&    4.7306\\
   45&   3.8023\\
   46&    4.7115\\
   48&    4.6959\\
   54&    6.3291\\
   56&    8.5134\\
   60&    7.2926\\
   65&    6.2505
   \end{tabular}
  \end{center}
     \vspace{-0.1cm}  
  \caption{Obtained optimal positions of the dampers and corresponding viscosities for Example \ref{exam:moderate}}
  \label{tab:Example2_pos}
 \end{table}

\section{Conclusion}\label{sec:conclusion}

We presented a new framework for damping optimization in mechanical systems 
described by a second-order differential equation, where the objective is to 
minimize the total average energy. The core idea is to penalize the total 
average energy criterion with a re-weighted $l_1$ norm, which promotes a sparse 
viscosity parameter and thereby allows the positions and viscosities of the 
dampers to be optimized simultaneously. After the viscosity optimization, we 
apply a pruning technique to estimate the optimal number of dampers, defined as 
the smallest number for which the relative increase in the objective function 
stays below a prescribed tolerance. We combined these ingredients into a single 
algorithm (SPARDO) and verified its efficiency and performance on grounded 
$n$-mass oscillators of small ($n = 20$) and moderate ($n = 100$) dimension. In 
both cases the estimated optimal number of dampers was $10\text{--}35\%$ of the 
system dimension, and where a brute-force comparison was feasible, SPARDO 
reproduced the same damper positions and objective values while being up to 
three orders of magnitude faster. Since the penalized step still optimizes over 
all $s$ variables, the method is currently practical for problems of small to 
moderate variable dimension; extending the framework to large-scale systems is 
a natural direction for future work.

\section*{Additional informations}
\subsection*{Funding}
Supported in part by the Croatian Science Foundation under the project “Optimization of parameter-dependent systems with applications”, Grant No. IP-2025-02-4862, and by the European Union-NextGenerationEU project, Grant No. 581-UNIOS-55, OpHoMat, and by the University of Osijek project MEOD, grant no. MZU-2026-15.
\subsection*{Data Sharing}
Data sharing not applicable to this article as no datasets were generated or analysed during the current study.
\subsection*{Conflict of interest}
There is no conflict of interest.
 
\bibliographystyle{plain}
\bibliography{csc, mor, software, References}

\end{document}